\documentclass[preprint,12pt]{elsarticle}
\usepackage{amsmath}

\usepackage{amssymb}

\begin{document}

\begin{frontmatter}


\ead{aaalikhanov@gmail.com}

\title{Boundary value problems for the diffusion equation  of the
variable order in differential and difference settings}


\author{Anatoly A. Alikhanov}

\address{Kabardino-Balkarian State University, ul. Chernyshevskogo 175,
 Nalchik,  360004,   Russia}

\begin{abstract}
Solutions of boundary value problems for a diffusion equation of
fractional and variable order in differential and difference
settings are studied. It is shown that the method of the energy
inequalities is applicable to obtaining a priori estimates for these
problems exactly as in the classical case. The credibility of the
obtained results is verified by performing numerical calculations
for a test problem.

\end{abstract}


\end{frontmatter}



\section{Introduction}

Fractional calculus is used for the description of a large class of
physical and chemical processes that occur in media with fractal
geometry as well as in the mathematical modeling of economic and
social-biological phenomena
\cite{Nakh:03,Podlub:99,Hilfer:00,Kilbas:06,Uchaikin:08}.    In
general, a medium in which a process proceeds is not homogenous,
moreover, its properties may vary in time. Mathematical models
containing equations with variable order derivatives provide a more
accurate and realistic description of processes proceeding in such
complex media (see e.g. \cite{Atanack:11,coimb:03,lorenzo:02}).

Therefore, the development of numerical and analytical methods of
the theory of fractional order differential equations is an actual
and important problem.

   Numerical methods for solving variable order fractional
  differential equations with various kinds of the  variable order
   fractional derivative have been proposed
\cite{Liu:12_1,Liu:12_2,Liu:10,Liu:09,Liu:09_2,Podlub:00}.

 The positivity
of the fractional derivative operator has been proved in
\cite{Nakh:03} and this result allows to obtain a priori estimates
for solutions of a large class of boundary value problems for the
equations containing fractional derivatives. The authors of the
paper \cite{ShkhTau:06} have obtained a priori estimate for the
solution of the Dirichlet boundary value problem of a fractional
order diffusion equation in terms of a fractional Riemann--Liouville
integral. The fractional diffusion equation with the regularized
fractional derivative has been studied, for example, in
\cite{Koch:90}. In the papers \cite{Main:96,MainGor:07,Pskh:09}, the
diffusion-wave equation with Caputo and Riemann--Liouville
fractional derivatives has been studied. The difference schemes for
boundary value problems for the fractional diffusion equation both
in one and multidimensional  cases have been studied
\cite{ShkhTau:06,Shkh:96,ShkhLaf:09}. A priori estimates for the
difference problems obtained in \cite{ShkhTau:06,Shkh:96,ShkhLaf:09}
by using the maximum principle imply the stability and convergence
of the considered difference schemes.

Using the energy inequality method, a priori estimates for the
solution of the Dirichlet and  Robin boundary value problems for the
diffusion-wave equation with Caputo fractional derivative have been
obtained \cite{Alikh:10}. More references on fractional order
differential equations, including the diffusion-wave equation, can
be found, for example, in \cite{Pskh:05}.

\section{Boundary value problems in differential setting}

\subsection{The Dirichlet boundary value problem}

In rectangle $\bar Q_T=\{(x,t):0\leq x\leq l, 0\leq t\leq T\}$ let
us study the boundary value problem

\begin{equation}\label{ur1}
\partial_{0t}^{\alpha(x)} u=\frac{\partial}{\partial
x}\left(k(x,t)\frac{\partial u}{\partial x}\right)-q(x,t)u+f(x,t),\,
0<x<l, 0<t\leq T,
\end{equation}

\begin{equation}
u(0,t)=0,\quad u(l,t)=0,\quad 0\leq t\leq T, \label{ur2}
\end{equation}

\begin{equation}
u(x,0)=u_0(x),\quad 0\leq x\leq l. \label{ur3}
\end{equation}

Where $0<c_1\leq k(x,t)\leq c_2$, $q(x,t)\geq 0$,
$\partial_{0t}^{\alpha(x)}
u(x,t)=\int_0^tu_{\tau}(x,\tau)(t-\tau)^{-\alpha(x)}d\tau/\Gamma(1-\alpha(x))$
 is a  Caputo fractional derivative of order $\alpha(x)$ ,
$0<\alpha(x)<1$ , for  all  $x\in (0,l)$, $\alpha(x)\in C(0,T)$
\cite{MainGor:07,Cap:69}.

Suppose further the existence of a solution  $u(x,t)\in C^{2,1}(\bar
Q_T)$ for the problem  (\ref{ur1})--(\ref{ur3}), where $C^{m,n}$ is
the class of functions, continuous together with their partial
derivatives of the order  $m$ with respect to  $x$ and order $n$
with respect to  $t$ on $\bar Q_T$.

 The existence of the solution for the initial boundary value
problem of a number of fractional order differential equation has
been proved \cite{luchko:10,luchko:11,Meersch:09}.

Let us prove the following:

{\bf Lemma 1.} For any functions $v(t)$ and $w(t)$ absolutely
continuous on $[0,T]$, one has the equality:

$$
v(t)\partial_{0t}^{\beta} w(t)+ w(t)\partial_{0t}^{\beta}
v(t)=\partial_{0t}^{\beta}(v(t)w(t))+
$$

\begin{equation}
+\frac{\beta}{\Gamma(1-\beta)}\int\limits_{0}^{t}\frac{d\xi}{(t-\xi)^{1-\beta}}
\int\limits_{0}^{\xi}\frac{v'(\eta)d\eta}{(t-\eta)^\beta}\int\limits_{0}^{\xi}\frac{w'(s)ds}{(t-s)^\beta},
 \label{ur4}
\end{equation}
where $0<\beta<1$.

 {\bf Proof.} Let us consider the difference
$$
v(t)\partial_{0t}^{\beta} w(t)+ w(t)\partial_{0t}^{\beta}
v(t)-\partial_{0t}^{\beta}(v(t) w(t))=
$$

$$
=\frac{1}{\Gamma(1-\beta)}\int\limits_{0}^{t}\frac{w'(s)(v(t)-v(s))+v'(s)(w(t)-w(s))}{(t-s)^{\beta}}ds=
$$

$$
=\frac{1}{\Gamma(1-\beta)}\int\limits_{0}^{t}\frac{v'(s)ds}{(t-s)^{\beta}}\int\limits_{s}^{t}w'(\xi)d\xi+
\frac{1}{\Gamma(1-\beta)}\int\limits_{0}^{t}\frac{w'(s)ds}{(t-s)^{\beta}}\int\limits_{s}^{t}v'(\xi)d\xi=
$$

$$
=\frac{1}{\Gamma(1-\beta)}\int\limits_{0}^{t}w'(\xi)d\xi\int\limits_{0}^{\xi}\frac{v'(s)ds}{(t-s)^{\beta}}
+\frac{1}{\Gamma(1-\beta)}\int\limits_{0}^{t}v'(\xi)d\xi\int\limits_{0}^{\xi}\frac{w'(s)ds}{(t-s)^{\beta}}=
$$

$$
=\frac{1}{\Gamma(1-\beta)}\int\limits_{0}^{t}(t-\xi)^\beta
\frac{\partial}{\partial\xi}\left(\int\limits_{0}^{\xi}\frac{v'(\eta)d\eta}{(t-\eta)^{\beta}}
\int\limits_{0}^{\xi}\frac{w'(s)ds}{(t-s)^{\beta}}\right)d\xi=
$$

$$
=\frac{\beta}{\Gamma(1-\beta)}\int\limits_{0}^{t}\frac{d\xi}{(t-\xi)^{1-\beta}}
\int\limits_{0}^{\xi}\frac{v'(\eta)d\eta}{(t-\eta)^\beta}\int\limits_{0}^{\xi}\frac{w'(s)ds}{(t-s)^\beta}.
$$
The proof of the lemma is complete.

If $v(t)=w(t)$ then from the Lemma 1 one has the following:

{\bf Corollary 1.} For any function $v(t)$ absolutely continuous on
$[0,T]$, the following equality takes place:

\begin{equation}
v(t)\partial_{0t}^{\beta}
v(t)=\frac{1}{2}\partial_{0t}^{\beta}v^2(t)+
\frac{\beta}{2\Gamma(1-\beta)}\int\limits_{0}^{t}\frac{d\xi}{(t-\xi)^{1-\beta}}
\left(\int\limits_{0}^{\xi}\frac{v'(\eta)d\eta}{(t-\eta)^\beta}\right)^2,
 \label{ur5}
\end{equation}
where $0<\beta<1$.

Let us use the following notation:
$\|u\|_0^2=\int\limits_{0}^{l}u^2(x,t)dx$,
$D_{0t}^{-\beta}u(x,t)=\int\limits_{0}^{t}(t-s)^{\beta-1}u(x,s)ds/\Gamma(\beta)$
-- fractional Riemann--Liouville integral of order  $\beta$.

{\bf Theorem 1.}   
If $k(x,t)\in C^{1,0}(\bar Q_T)$, $q(x,t)$, $f(x,t)\in C(\bar Q_T)$,
$k(x,t)\geq c_1>0$, $q(x,t)\geq 0$ everywhere on $\bar Q_T$, then
the solution  $u(x,t)$ of the problem (\ref{ur1})--(\ref{ur3})
satisfies the a priori estimate:

$$
\int\limits_{0}^{l}D_{0t}^{\alpha(x)-1}u^2(x,t)dx+c_1\int\limits_{0}^{t}\|u_x(x,s)\|_0^2ds\leq
$$

\begin{equation}
\leq \frac{l^2}{2c_1}\int\limits_{0}^{t}\|f(x,s)\|_0^2ds+
\int\limits_0^l\frac{t^{1-\alpha(x)}}{\Gamma(2-\alpha(x))}u_0^2(x)dx.
 \label{ur6}
\end{equation}

{\bf Proof.} Let us multiply equation (\ref{ur1}) by $u(x,t)$ and
integrate the resulting relation over $x$ from $0$ to $l$:

$$
\int\limits_0^lu(x,t)\partial_{0t}^{\alpha(x)}u(x,t)dx-
\int\limits_0^lu(x,t)(k(x,t)u_x(x,t))_xdx+
$$

\begin{equation}
+\int\limits_0^lq(x,t)u^2(x,t)dx=\int\limits_0^lu(x,t)f(x,t)dx.
 \label{ur7}
\end{equation}
Then transform the terms in identity (\ref{ur7}) as:

$$
- \int\limits_0^lu(x,t)(k(x,t)u_x(x,t))_xdx=
\int\limits_0^lk(x,t)u_x^2(x,t)dx\geq c_1\|u_x(x,t)\|_0^2,
$$

$$
\left|\int\limits_0^lu(x,t)f(x,t)dx\right|\leq
\varepsilon\|u(x,t)\|_0^2+\frac{1}{4\varepsilon}\|f(x,t)\|_0^2,
\quad \varepsilon>0,
$$
Using the equality (\ref{ur5}) one obtains
$$
\int\limits_0^lu(x,t)\partial_{0t}^{\alpha(x)}u(x,t)dx\geq
\frac{1}{2}\int\limits_0^l\partial_{0t}^{\alpha(x)}u^2(x,t)dx.
$$
Taking into account the above performed transformations, from the
identity  (\ref{ur7})  one arrives at the inequality

\begin{equation}
\frac{1}{2}\int\limits_0^l\partial_{0t}^{\alpha(x)}u^2(x,t)dx+c_1\|u_x(x,t)\|_0^2\leq
\varepsilon\|u(x,t)\|_0^2+\frac{1}{4\varepsilon}\|f(x,t)\|_0^2.
 \label{ur8}
\end{equation}

Using the inequality $\|u(x,t)\|_0^2\leq (l^2/2)\|u_x(x,t)\|_0^2$,
from the  inequality (\ref{ur8}) at $\varepsilon=c_1/l^2$ one
obtains
\begin{equation}
\int\limits_0^l\partial_{0t}^{\alpha(x)}u^2(x,t)dx+c_1\|u_x(x,t)\|_0^2
\leq \frac{l^2}{2c_1}\|f(x,t)\|_0^2.
 \label{ur911}
\end{equation}
Changing the variable $t$ by $s$ in the inequality (\ref{ur911}) and
integrating it over $s$ from $0$ to $t$, one obtains a priori
estimate (\ref{ur6}). The  uniqueness and the continuous dependence
of the solution of the problem (\ref{ur1})--(\ref{ur3}) on the input
data follow from the a priori estimate (\ref{ur6}).

The solution of the problem (\ref{ur1})--(\ref{ur3}) with
$\alpha(x)=\alpha$ ($\alpha=const$) satisfies the a priori
estimates:

\begin{equation}
D_{0t}^{\alpha-1}\|u(x,t)\|_0^2+c_1\int\limits_{0}^{t}\|u_x(x,s)\|_0^2ds
\leq\frac{l^2}{2c_1}\int\limits_{0}^{t}\|f(x,s)\|_0^2ds+
\frac{t^{1-\alpha}}{\Gamma(2-\alpha)}\|u_0(x)\|_0^2,
 \label{ur101}
\end{equation}

\begin{equation}
\|u(x,t)\|_0^2+D_{0t}^{-\alpha}\|u_x(x,t)\|_0^2
 \leq M\left(D_{0t}^{-\alpha}\|f(x,t)\|_0^2+\|u_0^2(x)\|_0^2\right).
 \label{ur1111}
\end{equation}

Inequality (\ref{ur101}) follows from  (\ref{ur6}), and the a priori
estimate  (\ref{ur1111}) follows from inequality  (\ref{ur8}) with
$\alpha(x)=\alpha$. Actually, applying the fractional integration
operator $D_{0t}^{-\alpha}$ to the both sides of inequality
(\ref{ur8}), one arrives at the estimate  (\ref{ur1111}), which
contains the constant $M=\max\{l^2/c_1 ,1\}/\min\{1,c_1\}$.

\subsection{The Robin boundary value problem.}

In the problem (\ref{ur1})--(\ref{ur3}) we replace  the boundary
conditions (\ref{ur2}) with
\begin{equation}
\left\{
\begin{array}{rcl}
k(0,t)u_x(0,t)=\beta_1(t)u(0,t)-\mu_1(t),     \\
-k(l,t)u_x(l,t)=\beta_2(t)u(l,t)-\mu_2(t).
\end{array}
\right. \label{ur10}
\end{equation}

In the rectangle $\bar Q_T$ we consider the Robin boundary value
problem (\ref{ur1}), (\ref{ur3}), (\ref{ur10}).

 {\bf Theorem 2.} If
$k(x,t)\in C^{1,0}(\bar Q_T)$, $q(x,t)$, $f(x,t)\in C(\bar Q_T)$,
$k(x,t)\geq c_1>0$, $q(x,t)\geq 0$ everywhere on $\bar Q_T$,
$\beta_i(t), \mu_i(t)\in C[0,T]$, $\beta_i(t)\geq \beta_0>0$, for
all $t\in [0,T]$, $i=1,2$, then the solution $u(x,t)$ of the problem
(\ref{ur1}), (\ref{ur3}), (\ref{ur10}) satisfies the a priori
estimate:

$$
\int\limits_{0}^{l}D_{0t}^{\alpha(x)-1}u^2(x,t)dx+
\gamma\left(\int\limits_{0}^{t}\left(\|u_x(x,s)\|_0^2+
u^2(0,s)+u^2(l,s)\right)ds\right)\leq
$$

\begin{equation}
\leq
\frac{\delta}{\gamma}\left(\int\limits_{0}^{t}\left(\|f(x,s)\|_0^2+\mu_1^2(s)+\mu_2^2(s)\right)ds\right)
+\int\limits_0^l\frac{t^{1-\alpha(x)}}{\Gamma(2-\alpha(x))}u_0^2(x)dx,
 \label{ur9}
\end{equation}

where $\gamma=\min\{c_1,\beta_0\}$, $\delta=\max\{1+l,l^2\}$.

{\bf Proof.} Multiply the equation (\ref{ur1}) by $u(x,t)$ and
integrate the resulting relation over $x$ from $0$ to $l$:
\begin{equation}
\int\limits_0^lu\partial_{0t}^\alpha
udx-\int\limits_0^lu(ku_x)_xdx+\int\limits_0^lqu^2dx=\int\limits_0^lufdx.
 \label{ur11}
\end{equation}

Then, transform the terms of the identity (\ref{ur11}):
$$
\int\limits_0^lu(x,t)\partial_{0t}^{\alpha(x)}u(x,t)dx\geq
\frac{1}{2}\int\limits_0^l\partial_{0t}^{\alpha(x)}u^2(x,t)dx.
$$
$$
-\int\limits_0^lu(ku_x)_xdx= \beta_1(t)u^2(0,t)+\beta_2(t)u^2(l,t)-
\mu_1(t)u(0,t)-\mu_2(t)u(l,t)+\int\limits_0^lku_x^2dx,
$$
$$
\left|\int\limits_0^lufdx\right|\leq
\varepsilon\|u\|_0^2+\frac{1}{4\varepsilon}\|f\|_0^2,\quad
\varepsilon>0.
$$
From (\ref{ur11}), taking into account the transformations
performed, one arrives at the inequality

$$
\frac{1}{2}\int\limits_0^l\partial_{0t}^{\alpha(x)}u^2(x,t)dx+c_1\|u_x(x,t)\|_0^2
+\beta_0u^2(0,t)+\beta_0u^2(l,t)\leq
$$

\begin{equation}
\leq \mu_1(t)u(0,t)+\mu_2(t)u(l,t)+
\varepsilon\|u\|_0^2+\frac{1}{4\varepsilon}\|f\|_0^2.
 \label{ur12}
\end{equation}

Using the inequalities $\mu_1(t)u(0,t)\leq\varepsilon
u^2(0,t)+(4\varepsilon)^{-1}\mu_1^2(t)$,
$\mu_2(t)u(l,t)\leq\varepsilon
u^2(l,t)+(4\varepsilon)^{-1}\mu_2^2(t)$, $\varepsilon>0$;
$\|u(x,t)\|_0^2\leq l^2\|u_x(x,t)\|_0^2+l(u^2(0,t)+u^2(l,t))$ with
$\varepsilon={\gamma}/({2\delta})$, from (\ref{ur12}) one has the
following inequality
$$
\int\limits_0^l\partial_{0t}^{\alpha(x)}u^2(x,t)dx+\gamma\left(\|u_x(x,t)\|_0^2
+u^2(0,t)+u^2(l,t)\right)\leq
$$
\begin{equation}
\leq
\frac{\delta}{\gamma}\left(\|f(x,t)\|_0^2+\mu_1^2(t)+\mu_2^2(t)\right).
\label{ur13}
\end{equation}

Changing variable $t$ by $s$ in inequality (\ref{ur13}) and
integrating it over $s$ from $0$ to $t$, one obtains the a priori
estimate (\ref{ur9}).

The uniqueness and the continuous dependence of the solution of
problem (\ref{ur1}), (\ref{ur3}), (\ref{ur10}) on the input data
follow from  the a priori estimate (\ref{ur9}).

The solution of the problem (\ref{ur1}), (\ref{ur3}), (\ref{ur10})
with $\alpha(x)=\alpha$ ($\alpha=const$) satisfies the following a
priori estimates:
$$
D_{0t}^{\alpha-1}\|u(x,t)\|_0^2+
\gamma\left(\int\limits_{0}^{t}\left(\|u_x(x,s)\|_0^2+
u^2(0,s)+u^2(l,s)\right)ds\right)\leq
$$
\begin{equation}
\leq
\frac{\delta}{\gamma}\left(\int\limits_{0}^{t}\left(\|f(x,s)\|_0^2+\mu_1^2(s)+\mu_2^2(s)\right)ds\right)
+ \frac{t^{1-\alpha}}{\Gamma(2-\alpha)}\|u_0(x)\|_0^2,
 \label{ur14}
\end{equation}
$$
\|u(x,t)\|_0^2+D_{0t}^{-\alpha}\|u_x(x,t)\|_0^2\leq
M(D_{0t}^{-\alpha}\|f(x,t)\|_0^2+
$$
\begin{equation}
+D_{0t}^{-\alpha}\mu_1^2(t)+D_{0t}^{-\alpha}\mu_2^2(t)+\|u_0^2(x)\|_0^2),
 \label{ur15}
\end{equation}
Inequality (\ref{ur14}) follows from  (\ref{ur9}), and the a priori
estimate  (\ref{ur15}) follows from inequality  (\ref{ur13}) with
$\alpha(x)=\alpha$. Actually, applying the fractional integration
operator $D_{0t}^{-\alpha}$ to the both sides of inequality
(\ref{ur13}), one arrives at the estimate  (\ref{ur15}), which
contains the constant $M=\max\{\delta/\gamma ,1\}/\min\{1,\gamma\}$.

\section{Boundary value problems in difference setting}

\subsection{The Dirichlet boundary value problem}

Suppose that a solution $u(x,t)\in C^{4,3}( Q_T)$ of the problem
(\ref{ur1})--(\ref{ur3}) exists, and the coefficients of the
equation (\ref{ur1}) and the functions $f(x,t)$, $u_0(x)$ satisfy
the smoothness conditions, required for the construction of
difference schemes with the order of approximation $O(\tau+h^2)$.

In the rectangle $\bar Q_T$ we introduce the grid
$\bar\omega_{h\tau}=\bar\omega_{h}\times\bar\omega_{\tau}$, where
$\bar\omega_{h}=\{x_i=ih, i=0,1,...,N, hN=l\}$,
$\bar\omega_{\tau}=\{t_j=j\tau, j=0,1,...,j_0, \tau j_0=T\}$.

To problem (\ref{ur1})--(\ref{ur3}), we assign the difference
scheme:

 \begin{equation}\label{ur2.4}
\Delta_{0t_j}^{\alpha_i} y=\Lambda( \sigma y^{j+1}+(1-\sigma)y^{j})
+\varphi,\quad i=1,...,N-1,\quad j=1,...,j_0-1,
\end{equation}

\begin{equation}
y(0,t)=0,\quad y(l,t)=0,\quad j=0,...,j_0, \label{ur2.5}
\end{equation}

\begin{equation}
y(x,0)=u_0(x),\quad i=0,...,N, \label{ur2.6}
\end{equation}
where $\Lambda y=(ay_{\bar x})_x-dy$, $a=k(x_{i-1/2},\bar t)$,
$d=q(x_i,\bar t)$, $\varphi=f(x_i,\bar t)$, $\bar t=t_{j+1/2}$,
$0\leq \sigma\leq 1$, $\Delta_{0t_j}^{\alpha_i}
y=\sum_{s=0}^j(t_{j-s+1}^{1-{\alpha_i}}-t_{j-s}^{1-{\alpha_i}})y_{
t}^s/\Gamma(2-{\alpha_i})$ -- a difference analogue of the Caputo
fractional derivative of order $\alpha_i$, $\alpha_i=\alpha(x_i)$
\cite{ShkhTau:06}.

According to \cite{ShkhTau:06,Samar:77} the difference scheme
(\ref{ur2.4})--(\ref{ur2.6}) has the order of approximation
$O(\tau+h^2)$.

{\bf Lemma 2.} For every function $y(t)$ defined on the grid
$\bar\omega_{\tau}$ one has the inequalities

\begin{equation}\label{ur2.8}
 y^{j+1}\Delta_{0t}^\alpha y \geq \frac{1}{2}\Delta_{0t}^\alpha (y^2)
 +\frac{\tau^\alpha \Gamma(2-\alpha)}{2}(\Delta_{0t}^\alpha y)^2,
\end{equation}

\begin{equation}\label{ur2.9}
 y^{j}\Delta_{0t}^\alpha y \geq \frac{1}{2}\Delta_{0t}^\alpha (y^2)
 -\frac{\tau^\alpha \Gamma(2-\alpha)}{2(2-2^{1-\alpha})}(\Delta_{0t}^\alpha
 y)^2.
\end{equation}

{\bf Proof.} Inequality (\ref{ur2.8}) is equivalent to the
inequality
$$
y^{j+1}\Delta_{0t}^\alpha y - \frac{1}{2}\Delta_{0t}^\alpha
(y^2)-\frac{\tau^\alpha \Gamma(2-\alpha)}{2}(\Delta_{0t}^\alpha
y)^2=
\frac{1}{\Gamma(2-\alpha)}y^{j+1}\sum\limits_{s=0}^j(t_{j-s+1}^{1-\alpha}-t_{j-s}^{1-\alpha})y_{
t}^s-
$$
$$
-\frac{1}{\Gamma(2-\alpha)}\sum\limits_{s=0}^j(t_{j-s+1}^{1-\alpha}-t_{j-s}^{1-\alpha})y_{
t}^s\frac{y^{s+1}+y^{s}}{2}- \frac{\tau^\alpha
\Gamma(2-\alpha)}{2}(\Delta_{0t}^\alpha y)^2=
$$
$$
=\frac{1}{\Gamma(2-\alpha)}\sum\limits_{s=0}^j(t_{j-s+1}^{1-\alpha}-t_{j-s}^{1-\alpha})y_{
t}^s(y^{j+1}-\frac{y^{s+1}+y^{s}}{2})- \frac{\tau^\alpha
\Gamma(2-\alpha)}{2}(\Delta_{0t}^\alpha y)^2=
$$
$$
=\frac{1}{\Gamma(2-\alpha)}\sum\limits_{s=0}^j(t_{j-s+1}^{1-\alpha}-t_{j-s}^{1-\alpha})y_{
t}^s(\frac{\tau}{2}y_{t}^s+\sum\limits_{k=s+1}^jy_{ t}^k\tau)-
\frac{\tau^\alpha \Gamma(2-\alpha)}{2}(\Delta_{0t}^\alpha y)^2=
$$
$$
=\frac{\tau}{2\Gamma(2-\alpha)}\sum\limits_{s=0}^j(t_{j-s+1}^{1-\alpha}-t_{j-s}^{1-\alpha})(y_{
t}^s)^2 +\frac{1}{\Gamma(2-\alpha)}\sum\limits_{k=1}^jy_{
t}^k\tau\sum\limits_{s=0}^{k-1}(t_{j-s+1}^{1-\alpha}-t_{j-s}^{1-\alpha})y_{
t}^s-
$$
\begin{equation}\label{ur2.10}
-\frac{\tau^\alpha \Gamma(2-\alpha)}{2}(\Delta_{0t}^\alpha
y)^2\geq0.
\end{equation}

Here we consider the sums to be equal to zero if the upper summation
limit is less than the lower one.

Let us introduce the following notation:
$\sum_{s=0}^k(t_{j-s+1}^{1-\alpha}-t_{j-s}^{1-\alpha})y_{
t}^s=v^{k+1}$, then $y_{
t}^0=(t_{j+1}^{1-\alpha}-t_{j}^{1-\alpha})^{-1}v^1$,\quad $y_{
t}^k=\tau(t_{j-k+1}^{1-\alpha}-t_{j-k}^{1-\alpha})^{-1}v_{
t}^k,\quad k=1,2,...,j$.\\
Taking into account the introduced notation, we rewrite the
inequality (\ref{ur2.10}) as
$$
\frac{\tau}{2\Gamma(2-\alpha)}(t_{j+1}^{1-\alpha}-t_{j}^{1-\alpha})^{-1}(v^1)^2
+\frac{\tau}{2\Gamma(2-\alpha)}\sum\limits_{k=1}^j\tau^2(t_{j-k+1}^{1-\alpha}-t_{j-k}^{1-\alpha})^{-1}(v_{
t}^k)^2+
$$
$$
+\frac{1}{\Gamma(2-\alpha)}\sum\limits_{k=1}^j\tau^2(t_{j-k+1}^{1-\alpha}-t_{j-k}^{1-\alpha})^{-1}v_{
t}^kv^{k}- \frac{\tau^\alpha}{2\Gamma(2-\alpha)}(v^{j+1})^2=
$$
$$
=\frac{\tau}{2\Gamma(2-\alpha)}(t_{j+1}^{1-\alpha}-t_{j}^{1-\alpha})^{-1}(v^1)^2
+\frac{1}{2\Gamma(2-\alpha)}\sum\limits_{k=1}^j\tau(t_{j-k+1}^{1-\alpha}-t_{j-k}^{1-\alpha})^{-1}((v^{k+1})^2-(v^{k})^2)-
$$
$$
-\frac{\tau^\alpha}{2\Gamma(2-\alpha)}(v^{j+1})^2=
\frac{1}{2\Gamma(2-\alpha)}\sum\limits_{k=0}^{j-1}
\tau((t_{j-k+1}^{1-\alpha}-t_{j-k}^{1-\alpha})^{-1}-
$$

\begin{equation}\label{ur2.11}
-(t_{j-k}^{1-\alpha}-t_{j-k-1}^{1-\alpha})^{-1})(v^{k+1})^2\geq 0.
\end{equation}

Obviously, inequality (\ref{ur2.11})  is valid since
$(t_{j-k+1}^{1-\alpha}-t_{j-k}^{1-\alpha})^{-1}-(t_{j-k}^{1-\alpha}-t_{j-k-1}^{1-\alpha})^{-1}>0,
k=0,1,...,j-1$.

 Let us prove now the inequality
(\ref{ur2.9}). Since $y^j=y^{j+1}-\tau y_t$, one obtains
$$
 y^{j}\Delta_{0t}^\alpha y - \frac{1}{2}\Delta_{0t}^\alpha (y^2)
 +\frac{\tau^\alpha \Gamma(2-\alpha)}{2(2-2^{1-\alpha})}(\Delta_{0t}^\alpha
 y)^2=
$$
$$
=y^{j+1}\Delta_{0t}^\alpha y - \frac{1}{2}\Delta_{0t}^\alpha (y^2)
 +\frac{\tau^\alpha \Gamma(2-\alpha)}{2(2-2^{1-\alpha})}(\Delta_{0t}^\alpha
 y)^2-\tau y_t\Delta_{0t}^\alpha y =
$$
$$
=\frac{\tau^\alpha(3-2^{1-\alpha})}{2\Gamma(2-\alpha)(2-2^{1-\alpha})}(v^{j+1})^2
-\frac{\tau^{1+\alpha}}{\Gamma(2-\alpha)}v_t^jv^{j+1} +
\frac{1}{2\Gamma(2-\alpha)}\sum\limits_{k=0}^{j-1}\tau
((t_{j-k+1}^{1-\alpha}-t_{j-k}^{1-\alpha})^{-1}-
$$
$$
-(t_{j-k}^{1-\alpha}-t_{j-k-1}^{1-\alpha})^{-1})(v^{k+1})^2
=\frac{\tau^\alpha(2^{1-\alpha}-1)}{2\Gamma(2-\alpha)(2-2^{1-\alpha})}(v^{j+1})^2
+\frac{\tau^\alpha}{\Gamma(2-\alpha)}v^{j+1}v^j+
$$
$$
+\frac{\tau^\alpha(2-2^{1-\alpha})}{2\Gamma(2-\alpha)(2^{1-\alpha}-1)}(v^{j})^2+
\frac{1}{2\Gamma(2-\alpha)}\sum\limits_{k=0}^{j-2}\tau
((t_{j-k+1}^{1-\alpha}-t_{j-k}^{1-\alpha})^{-1}-
$$
$$
-(t_{j-k}^{1-\alpha}-t_{j-k-1}^{1-\alpha})^{-1})(v^{k+1})^2\geq0.
$$
The proof of the lemma is complete.

{\bf Corollary 2.} For any function $y(t)$ defined on the grid
$\bar\omega_{\tau}$ one has the inequality
\begin{equation}\label{ur2.12}
 (\sigma y^{j+1}+(1-\sigma)y^{j})\Delta_{0t}^\alpha y \geq \frac{1}{2}\Delta_{0t}^\alpha
 (y^2)
 +\frac{\tau^\alpha \Gamma(2-\alpha)}{2(2-2^{1-\alpha})}\left((3-2^{1-\alpha})\sigma-1\right)(\Delta_{0t}^\alpha
 y)^2.
\end{equation}

{\bf Theorem 3.} The difference scheme (\ref{ur2.4})--(\ref{ur2.6})
at $\sigma\geq 1/(3-2^{1-\min(\alpha_i)})$ is absolutely stable and
its solution satisfies the following a priori estimate:

$$
\left(\frac{1}{\Gamma(2-\alpha_i)},\sum\limits_{j'=0}^{j}(t_{j-j'+1}^{1-\alpha_i}
-t_{j-j'}^{1-\alpha_i})(y^{j'+1})^2\right)+
$$

$$
+c_1\sum\limits_{j'=0}^{j}\|\sigma y_{\bar
x}^{j'+1}+(1-\sigma)y_{\bar x}^{j'}]|_0^2\tau\leq
$$

\begin{equation}\label{ur2.13.1}
\leq \frac{l^2}{2c_1}\sum\limits_{j'=0}^{j}\|\varphi^{j'}\|_0^2\tau+
\left(\frac{t_{j+1}^{1-\alpha_i}}{\Gamma(2-\alpha_i)},u_0^2(x_i)\right),
\end{equation}

where $(y,v)=\sum_{i=1}^{N-1}y_iv_ih$,
$(y,v]=\sum_{i=1}^{N}y_iv_ih$, $\|y\|_0^2=(y,y)$, $\|y]|_0^2=(y,y]$.

{\bf Proof.} Let us multiply scalarly equation (\ref{ur2.4}) by
$y^{(\sigma)}= \sigma y^{j+1}+(1-\sigma)y^{j}$:
\begin{equation}\label{ur2.7}
(\Delta_{0t_j}^{\alpha_i} y,y^{(\sigma)})-(\Lambda
y^{(\sigma)},y^{(\sigma)}) =(\varphi,y^{(\sigma)}).
\end{equation}
Let us transform the terms in identity (\ref{ur2.7}):
$$
(y^{(\sigma)},\Delta_{0t_j}^{\alpha_i} y)\geq
\frac{1}{2}(1,\Delta_{0t_j}^{\alpha_i} (y^2))+
$$
$$
 +(\frac{\tau^{\alpha_i}\Gamma(2-{\alpha_i})}{2(2-2^{1-{\alpha_i}})}
 \left((3-2^{1-{\alpha_i}})\sigma-1\right),(\Delta_{0t_j}^{\alpha_i}
 y)^2),
$$
$$
-(\Lambda y^{(\sigma)},y^{(\sigma)})=(a,(y_{\bar
x}^{(\sigma)})^2]+(d,(y^{(\sigma)})^2),
$$
$$
|(\varphi,y^{(\sigma)})|\leq
\varepsilon\|y^{(\sigma)}\|_0^2+\frac{1}{4\varepsilon}\|\varphi\|_0^2,
\quad\varepsilon>0.
$$
Taking into account the above-performed transformations, from
identity (\ref{ur2.7}) at $\sigma\geq 1/(3-2^{1-\min(\alpha_i)})$
one arrives at the inequality
\begin{equation}\label{ur2.13}
\frac{1}{2}(1,\Delta_{0t_j}^{\alpha_i} (y^2))+c_1\|y_{\bar
x}^{(\sigma)}]|_0^2\leq
\varepsilon\|y^{(\sigma)}\|_0^2+\frac{1}{4\varepsilon}\|\varphi\|_0^2.
\end{equation}

From (\ref{ur2.13}) at $\varepsilon=c_1/l^2$, using that
$\|y\|_0^2\leq (l^2/2)\|y_{\bar x}]|_0^2$, one obtains the
inequality

\begin{equation}\label{ur2.14}
(1,\Delta_{0t_j}^{\alpha_i} (y^2))+c_1\|y_{\bar
x}^{(\sigma)}]|_0^2\leq \frac{l^2}{2c_1}\|\varphi\|_0^2.
\end{equation}

Multiplying the inequality (\ref{ur2.14}) by $\tau$ and summing over
$j'$ from $0$ to $j$, one obtains the a priori estimate
(\ref{ur2.13.1}). The stability and convergence of the difference
scheme (\ref{ur2.4})--(\ref{ur2.6}) follow from the a priori
estimate (\ref{ur2.13.1}).

If $\alpha(x)=\alpha$ ($\alpha=const$) then the solution of the
problem (\ref{ur2.4})--(\ref{ur2.6}) satisfies the following a
priori estimate:
$$
\frac{1}{\Gamma(2-\alpha)}\sum\limits_{j'=0}^{j}(t_{j-j'+1}^{1-\alpha}-t_{j-j'}^{1-\alpha})\|y^{j'+1}\|_0^2
+c_1\sum\limits_{j'=0}^{j}\|\sigma y_{\bar
x}^{j'+1}+(1-\sigma)y_{\bar x}^{j'}]|_0^2\tau\leq
$$
\begin{equation}\label{ur2.15}
\leq \frac{l^2}{2c_1}\sum\limits_{j'=0}^{j}\|\varphi^{j'}\|_0^2\tau+
\frac{t_{j+1}^{1-\alpha}}{\Gamma(2-\alpha)}\|u_0(x_i)\|_0^2.
\end{equation}

Here the results are obtained for the homogeneous boundary
conditions  $u(0,t)=0$, $u(l,t)=0$. In the case of inhomogeneous
boundary conditions $u(0,t)=\mu_1(t)$, $u(l,t)=\mu_2(t)$ the
boundary conditions for the difference problem will have the
following form:
\begin{equation}
y(0,t)=\mu_1(t),\quad y(l,t)=\mu_2(t). \label{ur2.2.16}
\end{equation}

Convergence of the difference scheme (\ref{ur2.4}), (\ref{ur2.6}),
(\ref{ur2.2.16}) follows from the results obtained above. Actually,
let us introduce the notation $y=z+u$. Then the error  $z=y-u$ is a
solution of the following problem:

 \begin{equation}\label{ur2.2.17}
\Delta_{0t_j}^{\alpha_i} z=\Lambda( \sigma z^{j+1}+(1-\sigma)z^{j})
+\psi,\quad i=1,...,N-1,\quad j=1,...,j_0-1,
\end{equation}

\begin{equation}
z(0,t)=0,\quad z(l,t)=0,\quad j=0,...,j_0, \label{ur2.2.18}
\end{equation}

\begin{equation}
z(x,0)=0,\quad i=0,...,N, \label{ur2.2.19}
\end{equation}
where $\psi\equiv\Lambda( \sigma
u^{j+1}+(1-\sigma)u^{j})-\Delta_{0t_j}^{\alpha_i}
u+\varphi=O(\tau+h^2)$.

The solution of the problem (\ref{ur2.2.17})-(\ref{ur2.2.19})
satisfies the estimation (\ref{ur2.13.1}) so that the solution of
the difference scheme (\ref{ur2.4}), (\ref{ur2.6}), (\ref{ur2.2.16})
converges to the solution of the corresponding differential problem
with order $O(\tau+h^2)$.

\subsection{Numerical results}

 In this section, the following variable order time fractional
diffusion equation is considered:

\begin{equation}
\begin{cases}
\partial_{0t}^{\alpha(x)} u=\frac{\partial}{\partial
x}\left(k(x,t)\frac{\partial u}{\partial x}\right)-q(x,t)u+f(x,t),\,
0<x<1, 0<t\leq 1,\\
u(0,t)=\mu_1(t),\quad u(l,t)=\mu_2(t),\quad 0\leq t\leq 1,\\
u(x,0)=u_0(x),\quad 0\leq x\leq l,
\end{cases}
 \label{ur111}
\end{equation}
where $\alpha(x)=\frac{5+4\sin(6x)}{10}$,\quad
$k(x,t)=\frac{5+\cos(t)}{(3x^2+1)(3x^4+2x+1)+(x^3+x+1)(12x^3+2)}$,\quad
$q(x,t)=\frac{1+\sin(t)}{(x^3+x+1)(3x^4+2x+1)}$,\quad
$f(x,t)=(x^3+x+1)(3x^4+2x+1)(\frac{6t^{3-\alpha(x)}}{\Gamma(4-\alpha(x))}+\frac{6t^{2-\alpha(x)}}{\Gamma(3-\alpha(x))})+(1+\sin(t))(t^3+3t^2+1)$,
$\mu_1(t)=t^3+3t^2+1$,\quad $\mu_2(t)=18(t^3+3t^2+1)$,\quad
$u_0(x)=(x^3+x+1)(3x^4+2x+1)$.

The exact solution is $u(x,t)=(x^3+x+1)(3x^4+2x+1)(t^3+3t^2+1)$.

All the calculations are performed at $\sigma=
1/(3-2^{1-\min(\alpha_i)})$, where for the considered example
$\min(\alpha_i)=1/10$.

A comparison of the numerical solution and exact solution is
provided in {\bf Table 1}.

{\bf Table 2} shows that if $h=1/500$, then as the number of time of
our approximate scheme is decreased, a redaction in the maximum
error takes place, as expected and the convergence order of time is
$O(\tau)$, where the convergence order is given by the following
formula: Convergence order $
=\log_{\frac{\tau_1}{\tau_2}}\frac{e_1}{e_2}$.

{\bf Table 3} shows that when we take $h^2=\tau$, as the number as
spatial subintervals/time steps is decreased, a reduction in the
maximum error takes place, as expected the convergence order of the
approximate scheme is $O(h^2+\tau)$, where the convergence order is
given by the following formula: Convergence order $
=\log_{\frac{h_1}{h_2}}\frac{e_1}{e_2}$.

\begin{tabular}{lc}
{\bf Table 1}\\
The error, numerical solution and exact solution, when $t=1$, $h=1/10$, $\tau=1/100$.\\
\hline
\hspace{1mm} Space($x_i$) \hspace{15mm}{Numerical solution} \hspace{15mm}{Exact solution} \hspace{15mm}{Error}\\
\hline
\hspace{2mm} 0.0000 \hspace{19mm} 5.0000000 \hspace{30mm}  5.0000000 \hspace{22mm}    0.0000000  \\
\hspace{2mm} 0.1000 \hspace{19mm} 6.6068921 \hspace{30mm}  6.6076515 \hspace{22mm}    0.0007594  \\
\hspace{2mm} 0.2000 \hspace{19mm} 8.4813103 \hspace{30mm}  8.4849920 \hspace{22mm}    0.0036817  \\
\hspace{2mm} 0.3000 \hspace{17mm} 10.7677569 \hspace{28mm} 10.7772305 \hspace{22mm}   0.0094736  \\
\hspace{2mm} 0.4000 \hspace{17mm} 13.7191945 \hspace{28mm} 13.7381760 \hspace{22mm}   0.0189815  \\
\hspace{2mm} 0.5000 \hspace{17mm} 17.7407866 \hspace{28mm} 17.7734375 \hspace{22mm}   0.0326509  \\
\hspace{2mm} 0.6000 \hspace{17mm} 23.4561082 \hspace{28mm} 23.5063040 \hspace{22mm}   0.0501958  \\
\hspace{2mm} 0.7000 \hspace{17mm} 31.8041726 \hspace{28mm} 31.8738645 \hspace{22mm}   0.0696919  \\
\hspace{2mm} 0.8000 \hspace{17mm} 44.1761384 \hspace{28mm} 44.2609280 \hspace{22mm}   0.0847896  \\
\hspace{2mm} 0.9000 \hspace{17mm} 62.6016758 \hspace{28mm} 62.6793035 \hspace{22mm}   0.0776277  \\
\hspace{2mm} 1.0000 \hspace{17mm} 90.0000000 \hspace{28mm} 90.0000000 \hspace{22mm}   0.0000000  \\
 \hline
\end{tabular}

\vspace{5mm}

\begin{tabular}{lc}
{\bf Table 2}\\
Maximum error behavior versus time grid size reduction at $t=1$ when $h=1/500$.\\
\hline
\hspace{2mm} $\tau$ \hspace{39mm}{Maximum error} \hspace{40mm}{Convergence order} \\
\hline
\hspace{2mm} 1/256  \hspace{30mm} 0.0344960 \hspace{49mm}        \\
\hspace{2mm} 1/1024 \hspace{28mm} 0.0086690 \hspace{49mm}  0.996  \\
\hspace{2mm} 1/4096 \hspace{28mm} 0.0021738 \hspace{49mm}  0.998  \\
 \hline
\end{tabular}

\vspace{5mm}

\begin{tabular}{lc}
{\bf Table 3}\\
Maximum error behavior versus grid size reduction at $t=1$ when $h^2=\tau$.\\
\hline
\hspace{2mm} $h$ \hspace{39mm}{Maximum error} \hspace{40mm}{Convergence order} \\
\hline
\hspace{2mm} 1/40  \hspace{32mm} 0.0056275 \hspace{49mm}        \\
\hspace{2mm} 1/80 \hspace{32mm}  0.0014141 \hspace{49mm}  1.993  \\
\hspace{2mm} 1/160 \hspace{30mm} 0.0003542 \hspace{49mm}  1.997  \\
 \hline
\end{tabular}

\subsection{The Robin boundary value problem.}

To the differential problem (\ref{ur1}), (\ref{ur3}), (\ref{ur10})
we assign the following difference scheme:
\begin{equation}\label{ur4.1}
\Delta_{0t_j}^{\alpha_i} y=\Lambda( \sigma y^{j+1}+(1-\sigma)y^{j})
+\varphi, \quad i=0,...,N,  j=1,...,j_0,
\end{equation}

\begin{equation}
y(x,0)=u_0(x),\quad i=0,...,N, \label{ur4.2}
\end{equation}

where $\Lambda y=(a_1y_{x}-\tilde\beta_1 y)/(0.5h), i=0$, \quad
$\Lambda y=(ay_{\bar x})_x-dy, i=1,...,N-1$, \quad $\Lambda
y=(a_Ny_{\bar x}-\tilde\beta_2 y)/(0.5h), i=N$, \quad
$\varphi_0=(2\tilde\mu_1)/h$, \quad
$\varphi_N=(2\tilde\mu_2)/h$,\quad
$\tilde\beta_1=\beta_1+0.5hd_0$,\quad
$\tilde\beta_2=\beta_2+0.5hd_N$,\quad $\tilde\mu_1=\mu_1+0.5hf_0$,
\quad $\tilde\mu_2=\mu_2+0.5hf_N$.\\ The difference scheme
(\ref{ur4.1})--(\ref{ur4.2}) has the order of approximation
$O(\tau+h^2)$.

{\bf Theorem 4.} The difference scheme (\ref{ur4.1})--(\ref{ur4.2})
at $\sigma\geq 1/(3-2^{1-\min(\alpha_i)})$  is absolutely stable and
its solution satisfies the following a priori estimate:

$$
\left[\frac{1}{\Gamma(2-\alpha_i)},\sum\limits_{j'=0}^{j}
(t_{j-j'+1}^{1-\alpha_i}-t_{j-j'}^{1-\alpha_i})(y^{j'+1})^2\right]+
$$

$$
+ \gamma\sum\limits_{j'=0}^{j}\left(\|(y_{\bar
x}^{(\sigma)})^{j'+1}]|_0^2+((y_0^{(\sigma)})^{j'+1})^2+((y_N^{(\sigma)})^{j'+1})^2\right)\tau\leq
$$

$$
\leq \frac{\delta}{\gamma}\sum\limits_{j'=0}^{j}
\left((\tilde\mu_1^{j'+1/2})^2+(\tilde\mu_2^{j'+1/2})^2+\|\varphi^{j'}\|_0^2\right)\tau+
$$

\begin{equation}\label{ur4.2.2}+
\left[\frac{t_{j+1}^{1-\alpha_i}}{\Gamma(2-\alpha_i)},u_0^2(x_i)\right],
\end{equation}

where $\gamma=\min\{c_1,\beta_0\}$, $\delta=\max\{1+l,l^2\}$,
$[y,v]=\sum_{i=1}^{N-1}y_iv_ih+0.5y_0v_0h+0.5y_Nv_Nh$,
$|[y]|_0^2=[y,y]$, $(y^{(\sigma)})^{j'+1}= \sigma
y^{j'+1}+(1-\sigma)y^{j'}$.

{\bf Proof.} Let us multiply scalarly equation (\ref{ur4.1}) by
$y^{(\sigma)}= \sigma y^{j+1}+(1-\sigma)y^{j}$:
\begin{equation}\label{ur4.3}
[\Delta_{0t_j}^{\alpha_i} y,y^{(\sigma)}]-[\Lambda
y^{(\sigma)},y^{(\sigma)}] =[\varphi,y^{(\sigma)}],
\end{equation}
Let us transform the terms occurring in identity (\ref{ur4.3})
$$
[y^{(\sigma)},\Delta_{0t_j}^{\alpha_i} y]\geq
\frac{1}{2}[1,\Delta_{0t_j}^{\alpha_i} (y^2)]+
$$
$$
 +[\frac{\tau^{\alpha_i}\Gamma(2-{\alpha_i})}{2(2-2^{1-{\alpha_i}})}
 \left((3-2^{1-{\alpha_i}})\sigma-1\right),(\Delta_{0t_j}^{\alpha_i}
 y)^2],
$$
$$
-[\Lambda y^{(\sigma)},y^{(\sigma)}]=\tilde\beta_1
y_0^2+\tilde\beta_2 y_N^2+(a,(y_{\bar
x}^{(\sigma)})^2]+[d,(y^{(\sigma)})^2],
$$
$$
|[\varphi,y^{(\sigma)}]|\leq
\varepsilon\|y^{(\sigma)}\|_0^2+\tilde\mu_1y_0+\tilde\mu_2y_N+\frac{1}{4\varepsilon}\|\varphi\|_0^2,
\quad\varepsilon>0.
$$
Taking into account the above performed transformations, from
identity (\ref{ur4.3})  at $\sigma\geq 1/(3-2^{1-\min(\alpha_i)})$
one arrives at the inequality
$$
\frac{1}{2}[1,\Delta_{0t_j}^{\alpha_i} (y^2)]+c_1\|y_{\bar
x}^{(\sigma)}]|_0^2+\beta_0(y_0^2+y_N^2)\leq
$$
\begin{equation}\label{ur4.4}
\leq
\varepsilon(\|y^{(\sigma)}\|_0^2+y_0^2+y_N^2)+\frac{1}{4\varepsilon}(\tilde\mu_1^2+\tilde\mu_2^2+\|\varphi\|_0^2).
\end{equation}

From (\ref{ur4.3}) at $\varepsilon=\gamma/(2\delta)$, using that
$\|y\|_0^2\leq l^2\|y_{\bar x}]|_0^2+l(y_0^2+y_N^2)$, one has the
following inequality:
$$
[1,\Delta_{0t_j}^{\alpha_i} (y^2)]+\gamma(\|y_{\bar
x}^{(\sigma)}]|_0^2+y_0^2+y_N^2)\leq
$$
\begin{equation}\label{ur4.5}\leq
\frac{\delta}{\gamma}(\tilde\mu_1^2+\tilde\mu_2^2+\|\varphi\|_0^2).
\end{equation}
Multiplying inequality (\ref{ur4.5}) by $\tau$ and summing over $j'$
from $0$ to $j$, one obtains a priori estimate (\ref{ur4.2.2}). The
stability and convergence of the difference scheme
(\ref{ur4.1})--(\ref{ur4.2}) follow from the a priori estimate
(\ref{ur4.2.2}).

If $\alpha(x)=\alpha$ ($\alpha=const$), then for the solution of the
problem (\ref{ur4.1})--(\ref{ur4.2}) one has the following a priori
estimate:
$$
\frac{1}{\Gamma(2-\alpha)}\sum\limits_{j'=0}^{j}
(t_{j-j'+1}^{1-\alpha}-t_{j-j'}^{1-\alpha})|[y^{j'+1}]|_0^2+
$$
$$
+ \gamma\sum\limits_{j'=0}^{j}\left(\|(y_{\bar
x}^{(\sigma)})^{j'+1}]|_0^2+((y_0^{(\sigma)})^{j'+1})^2+((y_N^{(\sigma)})^{j'+1})^2\right)\tau\leq
$$
$$
\leq \frac{\delta}{\gamma}\sum\limits_{j'=0}^{j}
\left((\tilde\mu_1^{j'+1/2})^2+(\tilde\mu_2^{j'+1/2})^2+\|\varphi^{j'}\|_0^2\right)\tau+
$$
\begin{equation}\label{ur4.6}+
\frac{t_{j+1}^{1-\alpha}}{\Gamma(2-\alpha)}|[u_0(x_i)]|_0^2.
\end{equation}

\section{Conclusion}

The results obtained in the present paper allow to apply the method
of energy inequalities to finding a priori estimates for boundary
value problems for the fractional diffusion equation in differential
and difference settings exactly as in the classical case
($\alpha(x)=1$). It is interesting to note that the condition
$\sigma\geq 1/(3-2^{1-\alpha_i})$ at $\alpha(x)=1$ turns into the
well known condition $\sigma\geq 1/2$ of the absolute stability of
the difference schemes for the classical diffusion equation.

\section{Acknowledgements}

Dedicated to Prof. M. Kh. Shkhanukov, on the occasion of his 75-th
birthday.

\vskip 0.5cm

This work was supported by the Russian Foundation for Basic Research
(project 10-05-01150-a) and presented at the 4-th IFAC Workshop on
Fractional Differentiation and Its Applications, Badajoz, Spain,
October 18-20, 2010.






\begin{thebibliography}{00}

\bibitem{Nakh:03}
A.M. Nahushev, Fractional Calculus and its Application, {FIZMATLIT},
{Moscow}, {2003} (in Russian).


\bibitem{Podlub:99}
I. Podlubny, Fractional Differential Equations, Academic Press, San
Diego, 1999.

\bibitem{Hilfer:00}
R. Hilfer (Ed.), Applications of Fractional Calculus in Physics,
World Scientific, Singapore, 2000.

\bibitem{Kilbas:06}
A.A. Kilbas, H.M. Srivastava, J.J. Trujillo, Theory and Applications
of Fractional Differential Equation, Elsevier, Amsterdam, 2006.

\bibitem{Uchaikin:08}
V.V. Uchaikin, Method of Fractional Derivatives, Artishok,
Ul'janovsk, 2008 (in Russian).

\bibitem{Atanack:11}
T. M. Atanackovic, S. Pilipovic, {Hamilton's principle with variable
order fractional derivatives,} {Fract. Calc. Appl. Anal.} 14(1)
(2011) 94--109.

\bibitem{coimb:03}
C.F.M. Coimbra, Mechanics with variable-order differential
operators, Ann. Phys. 12 (11–12) (2003) 692-–703.

\bibitem{lorenzo:02}
C.F. Lorenzo, T.T. Hartley, Variable order and distributed order
fractional operators, Nonlinear Dynam. 29 (2002) 57–-98.

\bibitem{Liu:12_1}
S. Shen, F. Liu, J. Chen, I. Turner, V. Anh, {Numerical techniques
for the variable order time fractional diffusion equation,} {Appl.
Math. Comp.} {218} (2012) {10861--10870}.

\bibitem{Liu:12_2}
Chang-Ming Chen, F. Liu, V. Anh, I. Turner, {Numerical methods for
solving a two-dimensional variable-order anomalous subdiffusion
equation,} {Math. Comp.} {81} (2012) {345--366}.

\bibitem{Liu:10}
C. Chen, F. Liu, V. Anh,  I Turner, {Numerical schemes with high
spatial accuracy for a variable-order anomalous subdiffusion
equations,} {SIAM J. Scien. Comput.} 32(4) (2010) 1740--1760.

\bibitem{Liu:09}
P. Zhuang, F. Liu, V. Anh, I. Turner, {Numerical methods for the
variable-order fractional advection-diffusion equation with a
nonlinear source term,} {SIAM J. Numer. Anal.}  47(3) (2009)
1760--1781.

\bibitem{Liu:09_2}
R. Lin, F. Liu,  V. Anh, I. Turner, {Stability and convergence of
anewexplicitfinite-difference approximation for the variable-order
nonlinear fractional diffusion equation,} {Appl. Math. Comput.}  212
(2009) 435–-445.

\bibitem{Podlub:00}
I. Podlubny, Matrix approach to discrete fractional calculus, Fract.
Calc. Appl. Anal. 3(4) (2000)  359--386.


\bibitem{ShkhTau:06}
 M. Kh. Shkhanukov-Lafishev, F.I. Taukenova, {Difference methods for solving boundary value problems for fractional differential equations},
    {Comput. Math. Math. Phys.} 46(10) (2006)
    {1785--1795}.

\bibitem{Koch:90}
A. N. Kochubey, Diffusion of the fractional order, Differentsialnye
Uravneniya. 26(4) (1990) {660--770}.

 \bibitem{Main:96}
F. Mainardi, The fundamental solutions for the fractional
diffusion-wave equation, {Appl. Math. Lett.} {9(6)} (1996) 23--28.

\bibitem{MainGor:07}
F. Mainardi, R. Gorenflo, Time-fractional derivatives in relaxation
processes: a tutorial survey, Fract. Calc. Appl. Anal. {10(3)}
(2007) {269--308}.

\bibitem{Pskh:09}
A.V. Pskhu, The fundamental solution of a diffusion-wave equation of
fractional order, Izvestiya: Mathematics, 73(2) (2009) {351--392}
(in Russian).

\bibitem{Shkh:96}
M. Kh. Shkhanukov-Lafishev, {About convergence of difference schemes
for differential equations with fractional derivatives}, {Doklady
Akademii Nauk} {348(6)} (1996) {746--748} (in Russian).


 \bibitem{ShkhLaf:09}
M. Kh. Shkhanukov-Lafishev, M.M. Lafisheva, {Locally one-dimensional
difference schemes for the fractional order diffusion equation},
{Computational Mathematics and Mathematical Physics} {48(10)} (2009)
{1875--1884}.


\bibitem{Alikh:10}
A.A. Alikhanov, {A Priori Estimates for Solutions of Boundary Value
Problems for Fractional-Order Equations}, {Differ. Equ.} {46(5)}
(2010) {660--666}.


\bibitem{Pskh:05}
A.V. Pskhu, {Partial Differential Equations of the Fractional
Order}, {Nauka}, {Moscow}, {2005}.

 \bibitem{Cap:69}
M. Caputo, {Elasticita e Dissipazione}, {Zanichelli}, {Bologna},
{1969}.


\bibitem{luchko:10}
Y. Luchko, {Some uniqueness and existence results for the
initial-boundary-value problems for the generalized time-fractional
diffusion equation,} Comput. Math. Applic. 59 (2010)  1766--1772.

\bibitem{luchko:11}
Y. Luchko, {Initial-boundary-value problems for the generalized
multi-term time-fractional diffusion equation,} J. Math. Anal. Appl.
374 (2011) 538–-548.

\bibitem{Meersch:09}
M. Meerschaert, E. Nane, P. Vellaisamy, Fractional Cauchy problems
on bounded domains, Ann. Probab.  37 (2009) 979--1007.


\bibitem{Samar:77}
{A.A. Samarskiy,} {Theory of Difference Schemes}, {Nauka}, {Moscow},
{1977}.








 \end{thebibliography}


\end{document}